\documentclass[11pt]{amsart}

\usepackage{amscd,amssymb}
\usepackage{amsthm,amsmath,amssymb}

\usepackage[matrix,arrow]{xy}

\newtheorem{theorem}[equation]{Theorem}

\newtheorem{proposition}[equation]{Proposition}
\newtheorem{lemma}[equation]{Lemma}
\newtheorem{corollary}[equation]{Corollary}

\theoremstyle{definition}
\newtheorem{example}[equation]{Example}
\newtheorem{definition}[equation]{Definition}

\theoremstyle{remark}
\newtheorem{remark}[equation]{Remark}

\author{Ivan Cheltsov}

\title{Double cubics and double quartics}

\address{\begin{tabbing}\hspace*{28 em}\=\kill
Steklov Institute of Mathematics \>School of Mathematics\\
8 Gubkin street, Moscow 117966   \>University of Edinburgh\\
Russia                           \>Kings Buildings,  Mayfield Road\\
                                 \> Edinburgh EH9 3JZ, UK\\
\texttt{cheltsov@yahoo.com}      \>\\
                                 \>\texttt{I.Cheltsov@ed.ac.uk}
\end{tabbing}}

\begin{document}

\begin{abstract}
We study a double cover $\psi:X\to V\subset\mathbb{P}^{n}$
branched over a smooth divisor $R\subset V$ such that $R$ is cut
on $V$ by a hyper\-sur\-face of degree $2(n-\mathrm{deg}(V))$,
where $n\geqslant 8$ and $V$ is a smooth hypersurface of degree
$3$ or $4$. We prove that $X$ is non\-ra\-ti\-o\-nal and
birationally super\-ri\-gid.
\end{abstract}

\maketitle


\section{Introduction.}
\label{section:introduction}

Let $\psi:X\to V\subset\mathbb{P}^{n}$ be a double cover branched
over a smooth divisor $R\subset V$, where $n\geqslant 4$ and $V$
is a smooth hypersurface\footnote{All varieties are assumed to be
projective, normal, and defined over $\mathbb{C}$.}. Then
$\mathrm{rk}\,\mathrm{Pic}(X)=1$ (see \cite{Do82}) and
$$
-K_{X}\sim
\psi^{*}(\mathcal{O}_{\mathbb{P}^{n}}(d+r-1-n)\vert_{V}),
$$
where $d=\mathrm{deg}\,V$ and $r$ is a natural number such that
$R\sim \mathcal{O}_{\mathbb{P}^{n}}(2r)\vert_{V}$. Therefore $X$
is nonrational in the case when $d+r\geqslant n+1$. The variety
$X$ is rationally connected if $d+r\leqslant n$, because it is a
smooth Fano variety (see \cite{Ko96}). Moreover, the following
result is due to \cite{Pu00b}.

\begin{theorem}
\label{theorem:Pukhlikov}%
The variety $X$ is birationally superrigid\,\footnote{Namely, we
have $\mathrm{Bir}(X)=\mathrm{Aut}(X)$, and $X$ is not
bi\-ra\-ti\-o\-nal to the following varieties: a variety $Y$ such
that there is a morphism $\tau:Y\to Z$ whose general fiber has
negative Kodaira dimension and
$\mathrm{dim}(Y)\ne\mathrm{dim}(Z)\ne 0$; a Fano variety of Picard
rank $1$ having terminal $\mathbb{Q}$-factorial singularities that
is not biregular to $X$.} if it is general and $d+r=n\geqslant 5$.
\end{theorem}

In this paper we prove the following result.

\begin{theorem}
\label{theorem:main}%
The variety $X$ is birationally superrigid if $d+r=n\geqslant 8$
and $d=3$ or $4$.
\end{theorem}

One can use Theorem~\ref{theorem:main} to construct explicit
examples of nonrational Fano varieties.

\begin{example}
\label{example:double-quartic} The complete intersection
$$
\sum_{i=0}^{8}x_{i}^{4}=z^{2}-x_{0}^{4}x_{1}^{4}+x_{2}^{4}x_{3}^{4}+x_{4}^{4}x_{5}^{4}+x_{6}^{4}x_{7}^{4}=0\subset\mathbb{P}(1^{9},3)\cong\mathrm{Proj}(\mathbb{C}[x_{0},\ldots,x_{8},z])%
$$
is smooth. Hence, it is birationally superrigid and nonrational by
Theorem~\ref{theorem:main}.
\end{example}

In the case when $d+r=n\geqslant 4$ and $d=1$ or $2$ the
birational superrigidity of $X$ is proved in \cite{Is80b} and
\cite{Pu88a}. In the case when $d+r=n=4$ and $d=3$ the variety $X$
is not bi\-ra\-ti\-o\-nally superrigid, but it is nonrational
(see \cite{IsPu96}, \cite{Co00}). In the case when $d+r<n$ the
only known way to prove the nonrationality of $X$ is the method of
\S V in \cite{Ko96}, which implies the following result.

\begin{proposition}
\label{proposition:main}%
The variety $X$ is nonrational if it is very general, $n\geqslant
4$ and $r\geqslant {\frac{d+n+2}{2}}$.
\end{proposition}

The author would like to thank A.\,Corti, M.\,Grinenko,
V.\,Is\-kov\-skikh, J.\,Park, Yu.\,Pro\-kho\-rov and
V.\,Sho\-ku\-rov for useful and fruitful conversations.

\section{Preliminaries.}
\label{section:preliminaries}

Let $X$ be a variety and $B_X=\sum_{i=1}^{\epsilon}a_{i}B_{i}$ be
a boundary on $X$, where $a_{i}\in\mathbb{Q}$ and $B_{i}$ is
either a prime divisor on $X$ or a linear system on $X$ having no
base components. We say that $B_{X}$ is effective if every
$a_{i}\geqslant 0$, we say that $B_{X}$ is movable if every
$B_{i}$ is a linear system having no fixed
components\footnote{Every effective movable log pair can be
considered as a usual log pair (see \cite{KMM}).}. In the rest of
the section we we assume that all varieties are
$\mathbb{Q}$-factorial.

\begin{remark}
\label{remark:square-of-movable-boundary} We can consider
$B_{X}^{2}$ as an effective codimension-two cycle if $B_{X}$ is
movable.
\end{remark}

The notions such as discrepancies, terminality, cano\-ni\-city,
log terminality and log cano\-ni\-city can be defined for the log
pair $(X, B_{X})$ as for usual log pairs (see \cite{KMM}).

\begin{definition}
\label{definition:canonical-singularities} The log pair $(X,
B_{X})$ has canonical (terminal, respectively) singularities if
for every birational morphism $f:W\to X$ there is an equivalence
$$
K_{W}+B_{W}\sim_{\mathbb{Q}} f^{*}(K_{X}+B_{X})+\sum_{i=1}^{n}a(X,
B_{X}, E_{i})E_{i}
$$
such that every number $a(X, B_{X}, E_{i})$ is non-negative
(positive, respectively), where $B_{W}$ is a proper transform of
$B_{X}$ on $W$, and $E_{i}$ is an $f$-exceptional divisor. The
number $a(X, B_{X}, E_{i})$ is called the discrepancy of the log
pair $(X, B_X)$ in the divisor $E_i$.
\end{definition}

The application of Log Minimal Model Program (see \cite{KMM}) to
an effective movable log pair having canonical or terminal
singularities preserves its canonicity or terminality
respectively.

\begin{definition}
\label{definition:center-of-canonical-singularities} An
irreducible subvariety $Y\subset X$ is a center of canonical
singularities of the log pair $(X, B_{X})$ if there is a
birational morphism $f:W\to X$ and an $f$-ex\-cep\-tional divisor
$E$ such that $f(E)=Y$ and the inequality $a(X, B_{X}, E)\leqslant
0$ holds. The set of all centers of canonical
sin\-gu\-la\-ri\-ties of the log pair $(X, B_{X})$ is denoted as
$\mathbb{CS}(X, B_{X})$.
\end{definition}

In particular, the log pair $(X, B_{X})$ has terminal
singularities if and only if $\mathbb{CS}(X, B_{X})=\varnothing$.

\begin{remark}
\label{remark:canonical-reduction} Let $H$ be a general hyperplane
section of $X$. Then every component of $Z\cap H$ is contained in
the set $\mathbb{CS}(H, B_{X}\vert_{H})$ for every subvariety
$Z\subset X$ contained in $\mathbb{CS}(X, B_{X})$.
\end{remark}

\begin{remark}
\label{remark:canonical-centers-in-codimension-two} Let $Z\subset
X$ be a proper irreducible subvariety such that $X$ is smooth at
the generic point of $Z$. Suppose that $B_{X}$ is effective. Then
$Z\in \mathbb{CS}(X, B_{X})$ implies
$\mathrm{mult}_{Z}(B_{X})\geqslant 1$, but in the case
$\mathrm{codim}(Z\subset X)=2$ the inequality
$\mathrm{mult}_{Z}(B_{X})\geqslant 1$ implies $Z\in \mathbb{CS}(X,
B_{X})$.
\end{remark}

The following result is Lemma~3.18 in \cite{Che03}.

\begin{lemma}
\label{lemma:complete-intersection}%
Suppose that $X$ is a smooth complete intersection
$\cap_{i=1}^{k}G_{i}\subset\mathbb{P}^{n}$, and $B_{X}$ is
effective such that $B_{X}\sim_{\mathbb{Q}}rH$ for some
$r\in\mathbb{Q}$, where $G_{i}$ is a hypersurface in
$\mathbb{P}^{n}$, and $H$ is a hyperplane section of $X$. Then
$\mathrm{mult}_{Z}(B_{X})\leqslant r$ for every irreducible
subvariety $Z\subset X$ such that $\mathrm{dim}(Z)\geqslant k$.
\end{lemma}

The following result is well known (see \cite{Co95}, \cite{Co00}).

\begin{theorem}
\label{theorem:Nother-Fano} Let $X$ be a Fano variety of Picard
rank $1$ having terminal $\mathbb{Q}$-factorial singularities that
is not birationally superrigid. Then there is a linear system
$\mathcal{M}$ on the variety $X$ whose base locus has codimension
at least $2$ such that the singularities of the log pair $(X,
\mu\mathcal{M})$ are not canonical, where $\mu$ is a positive
rational number such that $K_{X}+\mu\mathcal{M}\sim_{\mathbb{Q}}
0$.
\end{theorem}

Let $f:V\to X$ be a birational morphism such that the union of
$\cup_{i=1}^{\epsilon} f^{-1}(B_{i})$ and all
$f$-ex\-ce\-pti\-onal divisors forms a divisor with simple normal
crossing. Then $f$ is called a log resolution of the log pair $(X,
B_{X})$, and the log pair $(V, B^{V})$ is called the log pull back
of $(X, B_{X})$ if
$$
B^{V}=f^{-1}(B_{X})-\sum_{i=1}^{n}a(X, B_{X}, E_{i})E_{i}
$$
such that $K_{V}+B^{V}\sim_{\mathbb{Q}}f^{*}(K_{X}+B_{X})$, where
$E_{i}$ is an $f$-exceptional divisor and $a(X, B_{X},
E_{i})\in\mathbb{Q}$.

\begin{definition}
\label{difinition:center-of-log-canonical-singularities} The log
canonical singularity subscheme $\mathcal{L}(X, B_{X})$ is the
subscheme associated to the ideal sheaf $\mathcal{I}(X,
B_{X})=f_{*}(\mathcal{O}_{V}(\lceil -B^{V}\rceil))$. A proper
irreducible subvariety $Y\subset X$ is called a center of log
canonical singularities of the log pair $(X, B_{X})$ if there is a
divisor $E\subset V$ that is contained in the effective part of
the support of $\lfloor B^{V}\rfloor$ and $f(E)=Y$. The set of all
centers of log canonical sin\-gu\-la\-ri\-ties of $(X, B_{X})$ is
denoted as $\mathbb{LCS}(X, B_{X})$, the set-theoretic union of
the elements of $\mathbb{LCS}(X, B_{X})$ is denoted as $LCS(X,
B_{X})$.
\end{definition}

In particular, we have $\mathrm{Supp}(\mathcal{L}(X,
B_{X}))=LCS(X, B_{X})$.

\begin{remark}
\label{remark:log-canonical-reduction} Let $H$ be a general
hyperplane section of $X$ and $Z\in \mathbb{LCS}(X, B_{X})$. Then
every component of the intersection $Z\cap H$ is contained in the
set $\mathbb{LCS}(H, B_{X}\vert_{H})$.
\end{remark}

The following result is Theorem~17.4 in \cite{Ko91}.

\begin{theorem}
\label{theorem:general-connectedness} Let $g:X\to Z$ be a
morphism. Then $LCS(X, B_{X})$ is connected in a neighborhood of
every fiber of the morphism $g\circ f$ if the following conditions
hold:
\begin{itemize}
\item the morphism $g$ has connected fibers;%
\item the divisor $-(K_{X}+B_{X})$ is $g$-nef and $g$-big;%
\item the inequality $\mathrm{codim}(g(B_{i})\subset Z)\geqslant 2$ holds if $a_{i}<0$;%
\end{itemize}
\end{theorem}

The following corollary of
Theorem~\ref{theorem:general-connectedness} is Theorem~17.6 in
\cite{Ko91}.

\begin{theorem}
\label{theorem:log-adjunction} Let $Z$ be an element of the set
$\mathbb{CS}(X, B_{X})$, and $H$ be an effective Cartier divisor
on the variety $X$. Suppose that the boundary $B_{X}$ is
effective, the varieties $X$ and $H$ are smooth in the generic
point of $Z$ and $Z\subset H\not\subset\mathrm{Supp}(B_{X})$. Then
$\mathbb{LCS}(H, B_{X}\vert_{H})\ne\varnothing$.
\end{theorem}

The following result is Theorem 3.1 in \cite{Co00}.

\begin{theorem}
\label{theorem:Corti} Suppose that $\mathrm{dim}(X)=2$, the
boundary $B_{X}$ is effective and movable, and there is a smooth
point $O\in X$ such that $O\in\mathbb{LCS}(X,
(1-a_{1})\Delta_{1}+(1-a_{2})\Delta_{2}+M_{X})$, where
$\Delta_{1}$ and $\Delta_{2}$ are smooth curves on $X$
intersecting normally at $O$, and $a_{1}$ and $a_{2}$ are
arbitrary non-negative rational numbers. Then we have
$$
\mathrm{mult}_{O}(B_{X}^{2})\geqslant\left\{\aligned
&4a_{1}a_{2}\ {\text {if}}\ a_{1}\leqslant 1\ {\text {or}}\ a_{2}\leqslant 1\\
&4(a_{1}+a_{2}-1)\ {\text {if}}\ a_{1}>1\ {\text {and}}\ a_{2}>1.\\
\endaligned
\right.
$$
\end{theorem}

\section{Main local inequality.}
\label{section:main-inequality}

Let $X$ be a variety, $O$ be a smooth point on $X$, $f:V\to X$ be
a blow up of the point $O$, $E$ be an exceptional divisor of $f$,
$B_X=\sum_{i=1}^{\epsilon}a_{i}\mathcal{B}_{i}$ be a movable
boundary on $X$, and $B_{V}=f^{-1}(B_{X})$, where $a_{i}$ is a
non-negative rational number and $\mathcal{B}_{i}$ is a linear
system on $X$ having no base components. Suppose that
$O\in\mathbb{CS}(X, B_{X})$, but the singularities of $(X, B_{X})$
is log terminal in some punctured neighborhood of the point $O$.
The following result is Corollary 3.5 in \cite{Co00}.

\begin{lemma}
\label{lemma:Corti} Suppose that $\mathrm{dim}(X)=3$ and
$\mathrm{mult}_{O}(B_{X})<2$. Then there is a line $L\subset
E\cong\mathbb{P}^{2}$ such that $L\in\mathbb{LCS}(V,
B_{V}+(\mathrm{mult}_{O}(B_{X})-1)E)$.
\end{lemma}

Suppose that $\mathrm{dim}(X)=4$ and $\mathrm{mult}_{O}(B_{X})<3$.
Then the proof of Lemma~\ref{lemma:Corti} and
Theorem~\ref{theorem:general-connectedness} implies the following
result.

\begin{proposition}
\label{proposition:Corti-dim-FOUR} One of the following
possibilities holds:
\begin{itemize}
\item there is a surface $S\subset E$ such that
$S\in\mathbb{LCS}(V,
B_{V}+(\mathrm{mult}_{O}(B_{X})-2)E)$;%
\item there is a line $L\subset E\cong\mathbb{P}^{3}$ such that
$L\in\mathbb{LCS}(V, B_{V}+(\mathrm{mult}_{O}(B_{X})-2)E)$.
\end{itemize}
\end{proposition}

Now suppose that the set $\mathbb{LCS}(V,
B_{V}+(\mathrm{mult}_{O}(B_{X})-2)E)$ does not contain surfaces
that are contained in the divisor $E$ and contains a line
$L\subset E\cong\mathbb{P}^{3}$. Let $g:W\to V$ be a blow up of
the variety $V$ in $L$, $F=g^{-1}(L)$, $\bar{E}=g^{-1}(E)$, and
$B_{W}=g^{-1}(B_{V})$. Then
$$
B^{W}=B_{W}+(\mathrm{mult}_{O}(B_{X})-3)\bar{E}+(\mathrm{mult}_{O}(B_{X})+\mathrm{mult}_{L}(B_{V})-5)F.
$$

\begin{proposition}
\label{proposition:CHELTSOV} One of the following possibilities
holds:
\begin{itemize}
\item the divisor $F$ is contained in $\mathbb{LCS}(W, B^{W}+\bar{E}+2F)$;%
\item there is a surface $Z\subset F$
such that $Z\in\mathbb{LCS}(W, B^{W}+\bar{E}+2F)$ and $g(Z)=L$.%
\end{itemize}
\end{proposition}

The following result is implied by
Proposition~\ref{proposition:CHELTSOV}.

\begin{theorem}
\label{theorem:technical} Let $Y$ be a variety,
$\mathrm{dim}(Y)=4$, $\mathcal{M}$ be a linear system on the
variety $Y$ having no base components, $S_{1}$ and $S_{2}$ be
sufficiently general divisors in $\mathcal{M}$, $P$ be a smooth
point on the variety $Y$ such that $P\in\mathbb{CS}(Y,
{\frac{1}{n}}\mathcal{M})$ for $n\in\mathbb{N}$, but the
singularities of $(Y, {\frac{1}{n}}\mathcal{M})$ are canonical in
some punctured neighborhood of the point $P$, $\pi:\hat{Y}\to Y$
be a blow up of $P$, and $\Pi$ be an exceptional divisor of $\pi$.
Then there is a line $C\subset\Pi\cong\mathbb{P}^{3}$ such that
the inequality
$$
\mathrm{mult}_{P}(S_{1}\cdot S_{2}\cdot \Delta)\geqslant 8n^{2}
$$
holds for any divisor $\Delta$ on $Y$ such that the following
conditions hold:
\begin{itemize}
\item the divisor $\Delta$ contains the point $P$ and $\Delta$ is smooth at $P$;%
\item the line $C\subset\Pi\cong\mathbb{P}^{3}$ is contained in the divisor $\pi^{-1}(\Delta)$;%
\item the divisor $\Delta$ does not contain subvarieties of dimension $2$ contained in $\mathrm{Bs}(\mathcal{M})$.%
\end{itemize}
\end{theorem}

\begin{proof}
Let $\Delta$ be a divisor on $Y$ such that $P\in\Delta$, the
divisor $\Delta$ is smooth at $P$, and $\Delta$ does not contain
any surface that is contained in the base locus of $\mathcal{M}$.
Then the base locus of the linear system
$\mathcal{M}\vert_{\Delta}$ has codimension $2$ in $\Delta$. In
particular, the intersection $S_{1}\cdot S_{2}\cdot \Delta$ is an
effective one-cycle. Let $\bar{S}_{1}=S_{1}\vert_{\Delta}$ and
$\bar{S}_{2}=S_{2}\vert_{\Delta}$. Then we must prove that the
inequality
\begin{equation}
\label{equation:inequality-8} \mathrm{mult}_{P}(\bar{S}_{1}\cdot \bar{S}_{2})\geqslant 8n^{2}%
\end{equation}
holds, perhaps, under certain additional conditions on $\Delta$.
Put $\bar{\mathcal{M}}=\mathcal{M}\vert_{\Delta}$. Then
$$
P\in\mathbb{LCS}(\Delta, {\frac{1}{n}}\bar{\mathcal{M}})
$$
by Theorem~\ref{theorem:log-adjunction}. Let
$\bar{\pi}:\hat{\Delta}\to\Delta$ be a blow up of $P$ and
$\bar{\Pi}=\bar{\pi}^{-1}(P)$. Then the diagram
\begin{equation}
\label{equation:commutative-diagram-of-blow-up} \xymatrix{
&\hat{\Delta}\ar@{->}[d]_{\bar{\pi}}\ar@{^{(}->}[rr]&&\hat{Y}\ar@{->}[d]^{\pi}&\\%
&\Delta\ar@{^{(}->}[rr]&&Y&} %
\end{equation}
is commutative, where $\hat{\Delta}$ is identified with
$\pi^{-1}(\Delta)\subset\hat{Y}$. We have
$\bar{\Pi}=\Pi\cap\hat{\Delta}$.

Let $\hat{\mathcal{M}}=\bar{\pi}^{-1}(\bar{\mathcal{M}})$. The
inequality~\ref{equation:inequality-8} is obvious if
$\mathrm{mult}_{P}(\bar{\mathcal{M}})\geqslant 3n$. Hence we may
assume that $\mathrm{mult}_{P}(\bar{\mathcal{M}})<3n$. Then
$$
\bar{\Pi}\not\in\mathbb{LCS}(\hat{\Delta}, {\frac{1}{n}}\hat{\mathcal{M}}+({\frac{1}{n}}\mathrm{mult}_{P}(\bar{\mathcal{M}})-2)\bar{\Pi}),%
$$
which implies the existence of a subvariety
$\Xi\subset\bar{\Pi}\cong\mathbb{P}^{2}$ such that $\Xi$ is a
center of log canonical singularities of $(\hat{\Delta},
{\frac{1}{n}}\hat{\mathcal{M}}+({\frac{1}{n}}\mathrm{mult}_{P}(\bar{\mathcal{M}})-2)\bar{\Pi})$.

Suppose that $\Xi$ is a curve. Put
$\hat{S}_{i}=\bar{\pi}^{-1}(S_{i})$. Then
$$
\mathrm{mult}_{P}(\bar{S}_{1}\cdot \bar{S}_{2})\geqslant\mathrm{mult}_{P}(\bar{\mathcal{M}})^{2}+\mathrm{mult}_{\Xi}(\hat{S}_{1}\cdot\hat{S}_{2}),%
$$
but we can apply Theorem~\ref{theorem:Corti} to the log pair
$(\hat{\Delta},
{\frac{1}{n}}\hat{\mathcal{M}}+({\frac{1}{n}}\mathrm{mult}_{P}(\bar{\mathcal{M}})-2)\bar{\Pi})$
in the generic point of the curve $\Xi$. The latter implies that
the inequality
$$
\mathrm{mult}_{\Xi}(\hat{S}_{1}\cdot\hat{S}_{2})\geqslant 4(3n^{2}-n\mathrm{mult}_{P}(\bar{\mathcal{M}}))%
$$
holds. Therefore we have
$$
\mathrm{mult}_{P}(\bar{S}_{1}\cdot\bar{S}_{2})\geqslant\mathrm{mult}_{P}(\bar{\mathcal{M}})^{2}+4(3n^{2}-n\mathrm{mult}_{P}(\bar{\mathcal{M}}))\geqslant 8n^{2},%
$$
which implies the inequality~\ref{equation:inequality-8}.

Suppose now that the subvariety $\Xi\subset\bar{\Pi}$ is a point.
In this case Proposition~\ref{proposition:Corti-dim-FOUR} implies
the existence of a line $C\subset\Pi\cong\mathbb{P}^{3}$ such that
$$
C\in\mathbb{LCS}(\hat{Y}, {\frac{1}{n}}\pi^{-1}(\mathcal{M})+(\mathrm{mult}_{P}(\mathcal{M})/n-2)\Pi)%
$$
and $\Xi=C\cap \hat{\Delta}$. The line $C\subset\Pi$ depends only
on the properties of the log pair $(Y, {\frac{1}{n}}\mathcal{M})$.

Suppose that initially we take $\Delta$ such that
$C\subset\pi^{-1}(\Delta)$. Then we can repeat all the previous
steps of our proof. Moreover, the geometrical meaning of
Proposition~\ref{proposition:CHELTSOV} is the following: the
condition $C\subset\hat{\Delta}=\pi^{-1}(\Delta)$ implies that
$$
C\in\mathbb{LCS}(\hat{\Delta}, {\frac{1}{n}}\hat{\mathcal{M}}+(\mathrm{mult}_{P}(\bar{\mathcal{M}})/n-2)\bar{\Pi})%
$$
in the case when the set $\mathbb{LCS}(\hat{\Delta},
{\frac{1}{n}}\hat{\mathcal{M}}+({\frac{1}{n}}\mathrm{mult}_{P}(\bar{\mathcal{M}})-2)\bar{\Pi})$
does not contain any other curve in $\bar{\Pi}$. Thus we can apply
the previous arguments to the divisor $\Delta$ such that
$C\subset\hat{\Delta}$ and obtain the proof of the
inequality~\ref{equation:inequality-8}.
\end{proof}

In the rest of the section we prove
Proposition~\ref{proposition:CHELTSOV}. We may assume that
$X\cong\mathbb{C}^{4}$. Let $H$ be a general hyperplane section of
$X$ such that $L\subset f^{-1}(H)$, $T=f^{-1}(H)$ and
$S=g^{-1}(T)$. Then
$$
K_{W}+B^{W}+\bar{E}+2F+S\sim_{\mathbb{Q}} (f\circ g)^{*}(K_{X}+B_{X}+H)%
$$
and
$$
B^{W}+\bar{E}+2F=B_{W}+(\mathrm{mult}_{O}(B_{X})-2)\bar{E}+(\mathrm{mult}_{O}(B_{X})+\mathrm{mult}_{L}(B_{V})-3)F,
$$
which implies that
$$
F\in\mathbb{LCS}(W, B^{W}+\bar{E}+2F)\iff \mathrm{mult}_{O}(B_{X})+\mathrm{mult}_{L}(B_{V})\geqslant 4%
$$
by
Definition~\ref{difinition:center-of-log-canonical-singularities}.
Thus we may assume that
$\mathrm{mult}_{O}(B_{X})+\mathrm{mult}_{L}(B_{V})<4$. We must
prove that there is a surface $Z\subset F$ such that
$Z\in\mathbb{LCS}(W, B^{W}+\bar{E}+2F)$ and $g(Z)=L$.

Now let $\bar{H}$ be a sufficiently general hyperplane section of
the variety $X$ passing through the point $O$,
$\bar{T}=f^{-1}(\bar{H})$ and $\bar{S}=g^{-1}(\bar{T})$. Then
$O\in\mathbb{LCS}(\bar{H}, B_{X}\vert_{{\bar H}})$ by
Theorem~\ref{theorem:log-adjunction} and
$$
K_{W}+B^{W}+\bar{E}+F+\bar{S}\sim_{\mathbb{Q}} (f\circ g)^{*}(K_{X}+B_{X}+H),%
$$
which implies that the log pair $(\bar{S},
(B^{W}+\bar{E}+F)\vert_{\bar{S}})$ is not log terminal. We can
apply Theorem~\ref{theorem:general-connectedness} to the morphism
$f\circ g:\bar{S}\to\bar{H}$. Therefore either the locus
$LCS(\bar{S}, (B^{W}+\bar{E}+F)\vert_{\bar{S}})$ consists of a
single isolated point in the fiber of the morphism
$g\vert_{F}:F\to L$ over the point $\bar{T}\cap L$ or it contains
a curve in the fiber of the morphism $g\vert_{F}:F\to L$ over the
point $\bar{T}\cap L$.

\begin{remark}
\label{remark:proof-of-main-inequality-I} Every element of the set
$\mathbb{LCS}(\bar{S}, (B^{W}+\bar{E}+F)\vert_{\bar{S}})$ that is
contained in the fiber of the $\mathbb{P}^{2}$-bundle
$g\vert_{F}:F\to L$ over the point $\bar{T}\cap L$ is an
intersection of $\bar{S}$ with some element of the set
$\mathbb{LCS}(W, B^{W}+\bar{E}+F)$ due to the generality in the
choice of $\bar{H}$.
\end{remark}

Therefore the generality of $\bar{H}$ implies that either
$\mathbb{LCS}(W, B^{W}+\bar{E}+F)$ contains a surface in the
divisor $F$ dominating the curve $L$ or the only center of log
canonical singularities of the log pair $(W, B^{W}+\bar{E}+F)$
that is contained in the divisor $F$ and dominates the curve $L$
is a section of the $\mathbb{P}^{2}$-bundle $g\vert_{F}:F\to L$.
On the other hand, we have
$$
\mathbb{LCS}(W, B^{W}+\bar{E}+F)\subseteq\mathbb{LCS}(W, B^{W}+\bar{E}+2F),%
$$
which implies that in order to prove
Proposition~\ref{proposition:CHELTSOV} we may assume that the
divisor $F$ contains a curve $C$ such that the following
conditions hold:
\begin{itemize}
\item the curve $C$ is a section of the $\mathbb{P}^{2}$-bundle $g\vert_{F}:F\to L$;%
\item the curve $C$ is the unique element of the set
$\mathbb{LCS}(W, B^{W}+\bar{E}+2F)$ that is contained in
the $g$-exceptional divisor $F$ and dominates the curve $L$;%
\item the curve $C$ is the unique element of the set
$\mathbb{LCS}(W, B^{W}+\bar{E}+F)$ that is contained in
the $g$-exceptional divisor $F$ and dominates the curve $L$.%
\end{itemize}

We have $O\in\mathbb{LCS}(H, M_{X}\vert_{H})$ by
Theorem~\ref{theorem:log-adjunction}, but $\mathbb{LCS}(S,
(B^{W}+\bar{E}+2F)\vert_{S})\ne\varnothing$, where $S$ is the
proper transform of $H$ on $W$. We can apply
Theorem~\ref{theorem:general-connectedness} to the log pair $(S,
(B^{W}+\bar{E}+2F)\vert_{S})$ and the birational morphism $f\circ
g\vert_{S}:S\to H$, which implies that one of the following holds:
\begin{itemize}
\item the locus $LCS(S, (B^{W}+\bar{E}+2F)\vert_{S})$ consists of a single point;%
\item the locus $LCS(S, (B^{W}+\bar{E}+2F)\vert_{S})$ contains a curve $C$.%
\end{itemize}

\begin{corollary}
\label{corollary:proof-of-main-inequality-II} Either $C\subset S$
or $S\cap C$ consists of a single point.
\end{corollary}

By construction we have $L\cong C\cong\mathbb{P}^{1}$ and
$$
F\cong\mathrm{Proj}(\mathcal{O}_{L}(-1)\oplus\mathcal{O}_{L}(1)\oplus\mathcal{O}_{L}(1))%
$$
and $S\vert_{F}\sim B+D$, where $B$ is the tautological line
bundle on $F$ and $D$ is a fiber of the natural projection
$g\vert_{F}:F\to L\cong\mathbb{P}^{1}$.

\begin{lemma}
\label{lemma:proof-of-main-inequality-III} The group
$H^{1}(\mathcal{O}_{W}(S-F))$ vanishes.
\end{lemma}

\begin{proof}
The intersection of the divisor $-g^{*}(E)-F$ with every curve
that is contained in the divisor $\bar{E}$ is non-negative and
$(-g^{*}(E)-F)\vert_{F}\sim B+D$. Hence $-4g^{*}(E)-4F$ is $h$-big
and $h$-nef, where $h=f\circ g$. However, we have
$X\cong\mathbb{C}^{4}$ and
$$
K_{W}-4g^{*}(E)-4F=S-F,
$$
which implies $H^{1}(\mathcal{O}_{W}(S-F))=0$ by the
Kawamata--Viehweg vanishing (see \cite{KMM}).
\end{proof}

Thus the restriction map
$$
H^{0}(\mathcal{O}_{W}(S))\to H^{0}(\mathcal{O}_{F}(S\vert_{F}))
$$
is surjective, but $\vert S\vert_{F}\vert$ has no base points (see
\S 2.8 in \cite{Re97}).

\begin{corollary}
\label{corollary:proof-of-main-inequality-IV} The curve $C$ is not
contained in $S$.
\end{corollary}

Let $\tau=g\vert_{F}$ and $\mathcal{I}_{C}$ be an ideal sheaf of
$C$ on $F$. Then $R^{1}\,\tau_{*}(B\otimes\mathcal{I}_{C})=0$ and
the map
$$
\pi:\mathcal{O}_{L}(-1)\oplus\mathcal{O}_{L}(1)\oplus\mathcal{O}_{L}(1)\to\mathcal{O}_{L}(k)%
$$
is surjective, where $k=B\cdot C$. The map $\pi$ is given by a an
element of the group
$$
H^{0}(\mathcal{O}_{L}(k+1))\oplus H^{0}(\mathcal{O}_{L}(k-1))\oplus H^{0}(\mathcal{O}_{L}(k-1)),%
$$
which implies $k\geqslant -1$.

\begin{lemma} \label{lemma:proof-of-main-inequality-V} The
equality $k=0$ is impossible.
\end{lemma}

\begin{proof}
Suppose $k=0$. Then the map $\pi$ is given by matrix $(ax+by, 0,
0)$, where $a$ and $b$ are complex numbers and $(x:y)$ are
homogeneous coordinates on $L\cong\mathbb{P}^{1}$. Thus the map
$\pi$ is not surjective over the point of $L$ at which $ax+by$
vanishes.
\end{proof}

Therefore the divisor $B$ can not have trivial intersection with
$C$. Hence the intersection of the divisor $S$ with the curve $C$
is either trivial or consists of more than one point, but we
already proved that $S\cap C$ consists of one point. The obtained
contradiction proves Proposition~\ref{proposition:CHELTSOV}.

The following result is a generalization of
Theorem~\ref{theorem:technical}.

\begin{theorem}
\label{theorem:technical-2} Let $Y$ be a variety of dimension
$r\geqslant 5$, $\mathcal{M}$ be a linear system on $Y$ having no
base components, $S_{1}$ and $S_{2}$ be general divisors in the
linear system $\mathcal{M}$, $P$ be a smooth point of the variety
$Y$ such that $P\in\mathbb{CS}(Y, {\frac{1}{n}}\mathcal{M})$ for
some natural number $n$, but the singularities of the log pair
$(Y, {\frac{1}{n}}\mathcal{M})$ are canonical in some punctured
neighborhood of $P$, $\pi:\hat{Y}\to Y$ be a blow up of the point
$P$, and $\Pi$ be a $\pi$-exceptional divisor. Then there is a
li\-ne\-ar subspace $C\subset\Pi\cong\mathbb{P}^{r-1}$ having
codimension $2$ such that $\mathrm{mult}_{P}(S_{1}\cdot S_{2}\cdot
\Delta)>8n^{2}$, where $\Delta$ is a divisor on $Y$ passing
through $P$ such that $\Delta$ is smooth at $P$, the divisor
$\pi^{-1}(\Delta)$ contains $C$, the divisor $\Delta$ does not
contain any subvarieties of $Y$ of codimension $2$ that are
contained in the base locus of $\mathcal{M}$.
\end{theorem}

\begin{proof}
We consider only the case $r=5$. Let $H_{1}, H_{2}, H_{3}$ be
general hyperplane sections of the variety $Y$ passing through
$P$. Put $\bar{Y}=\cap_{i=1}^{3}H_{i}$ and
$\bar{\mathcal{M}}=\mathcal{M}\vert_{\bar{Y}}$. Then $\bar{Y}$ is
a surface, which is smooth at $P$, and $P\in\mathbb{LCS}(\bar{Y},
{\frac{1}{n}}\bar{\mathcal{M}})$ by
Theorem~\ref{theorem:log-adjunction}. Let $\pi:\hat{Y}\to Y$ be a
blow up of $P$, $\Pi$ be an exceptional divisor of  $\pi$, and
$\hat{\mathcal{M}}=\pi^{-1}(\mathcal{M})$. Then the set
$$
\mathbb{LCS}(\hat{Y}, {\frac{1}{n}}\hat{\mathcal{M}}+(\mathrm{mult}_{P}(\mathcal{M})/n-2)\Pi)%
$$
contains a subvariety $Z\subset\Pi$ such that
$\mathrm{dim}(Z)\geqslant 2$.

In the case $\mathrm{dim}(Z)=4$ the claim is obvious. In the case
$\mathrm{dim}(Z)=3$ we can proceed as in the proof of
Theorem~\ref{theorem:technical} to prove that
$$
\mathrm{mult}_{P}(S_{1}\cdot S_{2}\cdot \Delta)>8n^{2}
$$
for any divisor $\Delta$ on $Y$ such that the divisor $\Delta$
contains the point $P$, the divisor $\Delta$ is smooth at the
point $P$, the divisor $\Delta$ does not contain any subvariety
$\Gamma\subset Y$ of codimension $2$ that is contained in the base
locus of the linear system $\mathcal{M}$.

It should be pointed out that in the cases when
$\mathrm{dim}(Z)\geqslant 3$ we do not need to fix any linear
subspace $C\subset\Pi$ of codimension $2$ such that
$\pi^{-1}(\Delta)$ contains $C$. The latter condition is vacuous
posteriori when $\mathrm{dim}(Z)\geqslant 3$.

Suppose that $\mathrm{dim}(Z)=2$. Then the surface $Z$ is a linear
subspace of $\Pi\cong\mathbb{P}^{4}$ having codimension $2$ by
Theorem~\ref{theorem:general-connectedness}. Moreover, the surface
$Z$ does not depend on the choice of our divisors $H_{1}$,
$H_{2}$, $H_{3}$, because it depends only on the properties of the
log pair $(Y, {\frac{1}{n}}\mathcal{M})$.

Put $C=Z$. Let $H$ be a sufficiently general hyperplane section of
$Y$ passing through the point $P$, and $\Delta$ be a divisor on
$Y$ such that $\Delta$ contains point $P$, the divisor $\Delta$ is
smooth at the point $P$, the divisor $\pi^{-1}(\Delta)$ contains
$C$, the divisor $\Delta$ does not contain any subvariety of $Y$
of codimension $2$ contained in the base locus of the linear
system $\mathcal{M}$. Then
$$
\mathrm{mult}_{P}(S_{1}\cdot S_{2}\cdot \Delta)>8n^{2}\iff \mathrm{mult}_{P}(S_{1}\vert_{H}\cdot S_{2}\vert_{H}\cdot \Delta\vert_{H})>8n^{2}%
$$
due to the generality of $H$. However, we have
$\mathrm{mult}_{P}(S_{1}\vert_{H}\cdot S_{2}\vert_{H}\cdot
\Delta\vert_{H})>8n^{2}$ by Theorem~\ref{theorem:technical},
because $P\in\mathbb{CS}(H, \mu\mathcal{M}\vert_{H})$ for some
positive rational number $\mu<1/n$ by
Theorem~\ref{theorem:log-adjunction}.
\end{proof}

\section{Birational superrigidity.}
\label{section:birational-superrigidity}

In this section we prove Theorem~\ref{theorem:main}. Let
$\psi:X\to V\subset\mathbb{P}^{n}$ be a double cover branched over
a smooth divisor $R\subset V$ such that $n\geqslant 7$. Then
$R\sim \mathcal{O}_{\mathbb{P}^{n}}(2r)\vert_{V}$ for some
$r\in\mathbb{N}$, and
$$
-K_{X}\sim
\psi^{*}(\mathcal{O}_{\mathbb{P}^{n}}(d+r-1-n)\vert_{V}),
$$
where $d=\mathrm{deg}\,V$. Suppose that $d+r=n$ and $d=3$ or $4$.
Then the group $\mathrm{Pic}(X)$ is generated by the divisor
$-K_{X}$, and $(-K_{X})^{2}=2d\leqslant 8$. Suppose that $X$ is
not birationally superrigid. Then
Theorem~\ref{theorem:Nother-Fano} implies the existence of a
linear system $\mathcal{M}$ whose base locus has codimension at
least $2$ and the singularities of the log pair $(X,
{\frac{1}{m}}\mathcal{M})$ are not canonical, where $m$ is a
natural number such that the equivalence $\mathcal{M}\sim -mK_{X}$
holds. Hence the set $\mathbb{CS}(X, {\frac{1}{m}}\mathcal{M})$
contains a proper irreducible subvariety $Z\subset X$ such that
$Z\in\mathbb{CS}(X, \mu\mathcal{M})$ for some rational $\mu<1/m$.

\begin{corollary}
\label{corollary:simple-inequality} For a general
$S\in\mathcal{M}$ the inequality $\mathrm{mult}_{Z}(S)>m$ holds.
\end{corollary}

A priori we have $\mathrm{dim}(Z)\leqslant \mathrm{dim}(X)-2=n-3$.
We may assume that $Z$ has maximal dimension among subvarieties of
$X$ such that the singularities of the log pair $(X,
{\frac{1}{m}}\mathcal{M})$ are not canonical in their generic
points.

\begin{lemma}
\label{lemma:exclution-of-points} The inequality
$\mathrm{dim}(Z)\ne 0$ holds.
\end{lemma}

\begin{proof}
Suppose that $Z$ is a point. Let $S_{1}$ and $S_{2}$ be
sufficiently general divisors in the linear system $\mathcal{M}$,
$f:U\to X$ be a blow up of $Z$, and $E$ be an $f$-exceptional
divisor. Then Theorem~\ref{theorem:technical-2} implies the
existence of a linear subspace $\Pi\subset E\cong\mathbb{P}^{n-2}$
of codimension $2$ such that
$$
\mathrm{mult}_{Z}(S_{1}\cdot S_{2}\cdot D)>8m^{2} %
$$
holds for any $D\in |-K_{X}|$ such that $\Pi\subset f^{-1}(D)$,
the divisor $D$ is smooth at $Z$, and $D$ does not contain any
subvariety of $X$ of codimension $2$ that is contained in the base
locus of $\mathcal{M}$.

Let $\mathcal{H}$ be a linear system of hyperplane sections of the
hypersurface $V$ such that $H\in \mathcal{H}$ if and only if
$\Pi\subset (\psi\circ f)^{-1}(H)$. Then there is a linear
subspace $\Sigma\subset\mathbb{P}^{n}$ of dimension $n-3$ such
that the divisors in the linear system $\mathcal{H}$ is cut on $V$
by the hyperplanes in $\mathbb{P}^{n}$ that contains the linear
subspace $\Sigma$. Hence the base locus of the linear system
$\mathcal{H}$ consists of the intersection $\Sigma\cap V$, but we
have $\Sigma\not\subset V$ by the Lefschetz theorem. In
particular, $\mathrm{dim}(\Sigma\cap V)=n-4$.

Let $H$ be a general divisor in $\mathcal{H}$ and
$D=\psi^{-1}(H)$. Then $\Pi\subset f^{-1}(D)$, and $D$ is smooth
at the point $Z$. Moreover, the divisor $D$ does not contain any
subvariety $\Gamma\subset X$ of codimension $2$ that is contained
in the base locus of $\mathcal{M}$, because otherwise
$\psi(\Gamma)\subset\Sigma\cap V$, but
$\mathrm{dim}(\psi(\Gamma))=n-3$ and $\mathrm{dim}(\Sigma\cap
V)=n-4$. Let $H_{1}, H_{2},\ldots, H_{k}$ be general divisors in
$|-K_{X}|$ passing through the point $Z$, where
$k=\mathrm{dim}(Z)-3$. Then we have
$$
2dm^{2}=H_{1}\cdot\cdots\cdot H_{k}\cdot S_{1}\cdot S_{2}\cdot D\geqslant \mathrm{mult}_{Z}(S_{1}\cdot S_{2}\cdot D)>8m^{2}, %
$$
which is a contradiction.
\end{proof}

\begin{lemma}
\label{lemma:exclution-of-big-codimension} The inequality
$\mathrm{dim}(Z)\geqslant \mathrm{dim}(X)-4$ holds.
\end{lemma}

\begin{proof}
Suppose that $\mathrm{dim}(Z)\leqslant\mathrm{dim}(X)-5$. Let
$H_{1}, H_{2},\ldots, H_{k}$ be sufficiently general hyperplane
sections of the hypersurface $V\subset\mathbb{P}^{n}$, where
$k=\mathrm{dim}(Z)>0$. Put
$$
\bar{V}=\cap_{i=1}^{k}H_{i},\ \bar{X}=\psi^{-1}(\bar{V}),\ \bar{\psi}=\psi\vert_{\bar{X}}:\bar{X}\to\bar{V},%
$$
and $\bar{\mathcal{M}}=\mathcal{M}\vert_{\bar{X}}$. Then $\bar{V}$
is a smooth hypersurface of degree $d$ in $\mathbb{P}^{n-k}$,
$\bar{\psi}$ is a double cover branched over a smooth divisor
$R\cap\bar{V}$, $\bar{\mathcal{M}}$ has no base components, and
$\bar{V}$ does not contains linear subspaces of $\mathbb{P}^{n-k}$
of dimension $n-k-3$ by the Lefschetz theorem. Let $P$ be any
point of the intersection $Z\cap\bar{X}$. Then
$P\in\mathbb{CS}(\bar{X}, {\frac{1}{m}}\bar{\mathcal{M}})$ and we
can repeat the proof of Lemma~\ref{lemma:exclution-of-points} to
get a contradiction.
\end{proof}

\begin{lemma}
\label{lemma:exclution-of-codimension-two-subvariety-I} The
inequality $\mathrm{dim}(Z)\ne \mathrm{dim}(X)-2$ holds.
\end{lemma}

\begin{proof}
Suppose that $\mathrm{dim}(Z)=\mathrm{dim}(X)-2$. Let $S_{1}$ and
$S_{2}$ be sufficiently general divisors in the linear system
$\mathcal{M}$, and $H_{1}, H_{2},\ldots, H_{n-3}$ be general
divisors in $|-K_{X}|$. Then
$$
2dm^{2}=H_{1}\cdot\cdots\cdot H_{n-3}\cdot S_{1}\cdot S_{2}\geqslant\mathrm{mult}_{Z}(S_{1})\mathrm{mult}_{Z}(S_{2})(-K_{X})^{n-3}\cdot Z>m^{2}(-K_{X})^{n-3}\cdot Z,%
$$
because $\mathrm{mult}_{Z}(\mathcal{M})>m$. Therefore
$(-K_{X})^{n-3}\cdot Z<2d$. On the other hand, we have
$$
(-K_{X})^{n-3}\cdot Z=\left\{\aligned
&\mathrm{deg}(\psi(Z)\subset\mathbb{P}^{n})\ \text{when}\ \psi\vert_{Z}\ \text{is birational},\\
&2\mathrm{deg}(\psi(Z)\subset\mathbb{P}^{n})\ \text{when}\ \psi\vert_{Z}\ \text{is not birational}.\\
\endaligned
\right.
$$

The Lefschetz theorem implies that $\mathrm{deg}(\psi(Z))$ is a
multiple of $d$. Therefore $\psi\vert_{Z}$ is a birational
morphism and $\mathrm{deg}(\psi(Z))=d$. Hence either $\psi(Z)$ is
contained in $R$, or the scheme-theoretic intersection
$\psi(Z)\cap R$ is singular in every point. However, we can apply
the Lefschetz theorem to the smooth complete intersection
$R\subset\mathbb{P}^{n}$, which gives a contradiction.
\end{proof}

\begin{lemma}
\label{lemma:exclution-of-codimension-three-subvariety-I} The
inequality $\mathrm{dim}(Z)\leqslant\mathrm{dim}(X)-5$ holds.
\end{lemma}

\begin{proof}
Suppose that $\mathrm{dim}(Z)\geqslant\mathrm{dim}(X)-4\geqslant
3$. Let $S$ be a sufficiently general divisor in the linear system
$\mathcal{M}$, $\hat{S}=\psi(S\cap R)$ and $\hat{Z}=\psi(Z\cap
R)$. Then $\hat{S}$ is a divisor on the complete intersection
$R\subset\mathbb{P}^{n}$ such that
$\mathrm{mult}_{\hat{Z}}(\hat{S})>m$ and
$\hat{S}\sim\mathcal{O}_{\mathbb{P}^{n}}(m)\vert_{R}$, because $R$
is a ramification divisor of $\psi$. Hence, the inequality
$\mathrm{dim}(\hat{Z})\geqslant 2$ is impossible by
Lemma~\ref{lemma:complete-intersection}.
\end{proof}

Therefore Theorem~\ref{theorem:main} is proved.

\section{Reduction into characteristic $2$.}
\label{section:double-covers}

In this section we prove Proposition~\ref{proposition:main}. The
following result is Theorem~5.12 in \S V of \cite{Ko96}.

\begin{theorem}
\label{theorem:Matsusaka} Let $f:X\to S$ be a proper and flat
morphism having irreducible and reduced fibers, $g:Z\to T$ be a
proper and flat morphism having reduced fibers, where $S$ is
irreducible scheme, and $T$ is a spectrum of discrete valuation
ring with closed point $O$. Suppose that some component of the
fiber $g^{-1}(O)$ is not geometrically ruled and the generic fiber
of $g$ is birational to a fiber of the morphism $f$. Then there
are countably many closed subvarieties $S_{i}\subset S$ such that
for any closed point $s\in S$ the fiber $f^{-1}(s)$ is
geometrically ruled $\iff$ $s\in\cup S_{i}$.
\end{theorem}

Let $Y$ be a scheme, $L$ be a line bundle on the scheme $Y$, and
$s$ be a global section of the line bundle $L^{k}$ for some
$k\in\mathbb{N}$. Let us construct a $k:1$ cover $Y_{s,\,L}^{k}$
of $Y$ ramified along the zeroes of the section $s$ as follows:
\begin{itemize}
\item let $U$ be a total space of $L$ with a natural projection $\pi:U\to Y$; %
\item we have $\pi_{*}(\mathcal{O}_{U})=\oplus_{i\geqslant 0}L^{-i}$ and $\pi_{*}(\pi^{*}(L))=\oplus_{i\geqslant -1}L^{-i}$;%
\item there is a canonical section $y$ of $\pi^{*}(L)$ that corresponds to $1\in H^{0}(\mathcal{O}_{Y})$;%
\item both $y$ and $s$ can be viewed as a section of $\pi^{*}(L^{k})$ since $\pi_{*}(\pi^{*}(L^{k}))=\oplus_{i\geqslant -k}L^{-i}$;%
\item let $y^k=s$ be an equation of $Y_{s,\,L}^{k}$ in $U$;%
\item there is a natural projection $\pi\vert_{Y_{s,\,L}^{k}}:Y_{s,\,L}^{k}\to Y$;%
\item the morphism $\pi\vert_{Y_{s,\,L}^{k}}$ is a $k:1$ cover ramified along the zeroes of the section $s$.%
\end{itemize}

\begin{example}
\label{example:special-cover} Let $Y=\mathbb{P}^n$ considered as a
scheme over $\mathbb{Z}$, $L=\mathcal{O}_{\mathbb{P}^n}(r)$ for
some $r\in\mathbb{N}$, and $s$ be a global section of
$\mathcal{O}_{\mathbb{P}^n}(2r)$. Consider the weighted projective
space
$$
\mathbb{P}(1,\ldots,1,r)=\text{Proj}(\mathbb{Z}[x_0,\ldots,x_n, y])%
$$
where $\mathrm{wt}(y)=r$ and $\mathrm{wt}(x_i)=1$. Then
$Y^2_{s,\,L}\cong V(y^2-s)\subset\mathbb{P}(1,\ldots,1,r)$.%
\end{example}

The following result is Theorem~5.11 in \S V of \cite{Ko96}.

\begin{theorem}
\label{theorem:Janos} Let $Y$ be a smooth projective variety over
an algebraically closed field of characteristic $p$, $L$ be a line
bundle on $Y$, $s$ be a general global section of $L^{p}$ such
that $\mathrm{dim}(Y)\geqslant 3$, the divisor $L^{p}\otimes
K_{V}$ is ample and the restriction map $H^{0}(Y,
L^{p})\to(\mathcal{O}_{Y}\slash m_{x}^{4})\otimes L^{p}$ is
surjective for every point $x\in Y$. Then $Y_{s,\,L}^{p}$ is not
separably uniruled.
\end{theorem}

Let $Y$ be a smooth hypersurface in $\mathbb{P}^{n}$ of degree $d$
defined over an algebraically closed field of characteristic $2$.
Let $L=\mathcal{O}_{\mathbb{P}^{n}}(r)\vert_{V}$ for some
$r\in\mathbb{N}$ and $s$ be a sufficiently general global section
of the line bundle $\mathcal{O}_{\mathbb{P}^{n}}(2r)\vert_{V}$.
Then $Y_{s,\,L}^{2}$ is not ruled if $r\geqslant
{\frac{d+n+2}{2}}$ and $n\geqslant 4$ by
Theorem~\ref{theorem:Janos}, therefore,
Theorem~\ref{theorem:Matsusaka} implies
Proposition~\ref{proposition:main}.

\end{document}